\documentclass[12pt,preprint,spbasic]{aastex}
\usepackage{spr-astr-addons}
\usepackage{aas_macros}
\usepackage{natbib}
\usepackage{url}

\newcommand{\abc}{{\Bigl(1+\frac{5}{2}A_2\Bigr) }}
\newcommand{\aac}{{(1-\frac{A_2}{2}) }}
\newcommand{\amc}{{(1-\mu) }}
\newcommand{\adc}{{\frac{\delta ^2}{2} }}

\newcommand{\zd}{{\displaystyle}}
\newcommand{\zabs}{{\Bigl(1+\frac{5A_2}{2{r^2_2}_i}\Bigr)}}
\newcommand{\zws}{{\frac{W_1}{{r^2_1}_i}}}

\shorttitle{The Effect of  Radiation Pressure ---- Restricted Three Body Problem}

\shortauthors{Badam Singh Kushvah}

\begin{document}
\title{The Effect of  Radiation Pressure on the Equilibrium Points in the Generalized Photogravitational Restricted Three Body Problem}

\author{Badam Singh Kushvah\altaffilmark{1}}
\affil{Gwalior Engineering College  Airport  Road, Maharajpura,
Gwalior (M.P.)-474015,INDIA}\email{bskush@gmail.com}

 \altaffiltext{1}{Present
address for correspondence : C/O Er. Prashant Anand Kushwaha,
139/I, Anupam Nagar Ext-2, Opp. Jiwaji University, City Center,
Gwalior-474011} 

%

\begin{abstract}
The existence of equilibrium points and the effect of radiation pressure  have been discussed numerically. The problem is generalized  by considering bigger primary as a source of radiation and small primary as an oblate spheroid. We have also discussed  the Poynting-Robertson(P-R) effect which is caused due to radiation pressure. It is found that the collinear points $L_1,L_2,L_3$ deviate from the axis joining the two primaries, while the triangular points $L_4,L_5$  are not symmetrical due to radiation pressure. We have seen that $L_1,L_2,L_3$ are linearly unstable while  $L_4,L_5$ are conditionally stable in the sense of Lyapunov when P-R effect is not considered. We have  found that  the effect of radiation pressure reduces the  linear stability zones while P-R effect induces an instability in the sense of Lyapunov.
\end{abstract}


\keywords{Radiation Pressure: Equilibrium Points: Generalized
Photogravitational:RTBP:Linear Stability: P-R effect}

%
%
%

\section{Introduction}
\label{intro}
Three-Body problem is a continuous source of study, since
the discovery of its non-integrability due to \citet*{Poincare1892}. Many of the best minds in Mathematics and Physics worked on this problem in the last century. Regular and chaotic motions have been widely investigated with any
kind of tools, from analytical results to numerical explorations. The
restricted three body model(RTBP) is insoluble, although specific
solutions exist, like the ones in which the spacecraft is positioned
at one of the Lagrangian points. The five singular points of  Jacobian function are called equilibrium points, Lagrangian points and most frequently are referred to as the equidistant(triangular) and collinear(straight line) solutions. There are many periodic orbits in the restricted three body problem. One of  the most famous,
discovered by Lagrange, is formed by an equilateral triangle. In
\lq\lq Earth-Moon-Space Station" model, if the Moon is not too
massive,  this orbit is thought to be stable. If the space station
is pushed a bit to one side (in position or velocity ), it's
supposed to make small oscillations around  this orbit. The two
kinds of triangular points are called $L_4$ and $L_5$ points. We
study the motion of three finite bodies in the three body problem.
The problem is restricted in the sense that one of the masses is
taken to be so small that the gravitational effect on the other
masses by third mass is negligible. The smaller body is known as
infinitesimal mass and remaining two  finite massive bodies as
primaries. The classical restricted three body problem is generalized to include the force of radiation pressure, oblateness effect and Poynting-Robertson(P-R) effect. 
The solar radiation pressure force $F_p$ is exactly opposite to the gravitational attraction force $F_g$ and change with the distance by the same law it is possible to consider that the result of action of this force will lead to reducing the effective mass of the  Sun or particle. It is acceptable to speak about a reduced mass of the particle as the effect of reducing its mass depends on the properties of the particle itself.

\citet*{Chernikov1970} and \citet*{Schuerman1980}  discussed the
position as well as the stability of the Lagrangian equilibrium
points when radiation pressure, P-R drag force are included.
\citet*{Murray1994} systematically discussed the dynamical effect of
general drag in the planar circular restricted three body problem.
\citet*{KushvahBR2006} examined the linear stability of triangular
equilibrium points in the generalized photogravitational restricted
three body problem with Poynting-Robertson drag, $L_4$ and $L_5$
points became unstable due to P-R drag which is very remarkable and
important, where as they are linearly stable in classical problem when $0<\mu<\mu_{Routh}=0.03852$. \citet*{KushvahetalHr2007,KushvahetalNON2007,KushvahetalNr2007} examined normalization of Hamiltonian they have also studied the nonlinear stability of triangular equilibrium points in the generalized photogravitational restricted three body problem with Poynting-Robertson drag, they have found  that the triangular points are stable in the nonlinear sense except three critical mass ratios at which KAM theorem fails.
\section{Equations of Motion}
\label{sec:eqmot} 
\cite{Poynting1903} has stated that the particle such as small
meteors or cosmic dust  are comparably affected by gravitational and
light  radiation force, as they approach luminous celestial bodies.
He also suggested that infinitesimal body in solar orbit suffers a
gradual loss of angular momentum and ultimately spiral into the Sun.
In a system of coordinates where the Sun is at rest, radiation 
scattered by infinitesimal mass in the direction of motion suffers a
blue shift and in the opposite direction it is red shifted. This
gives rise to net drag force which opposes the direction of motion. 
The proper relativistic treatment of this problem  was formulated by
\cite{Robertson1937} who showed that to first order in
${\vec{\frac{V}{c}}}$ the radiation pressure force is given by
\begin{equation}\vec{F}=\displaystyle{{F_p}\biggl\{ \frac{\vec{R}}{R} - \frac{\vec{V}.\vec{R}\vec{R}}{cR^{2}}- \frac{\vec{V}}{c} \biggr\}}\label{eq:F}\end{equation}
Where $F_{p}$=$\frac{3Lm}{16{\pi}R^{2}\rho{sc}}$   denotes the
measure of the radiation pressure force, $\vec{R}$  the position
vector of $P$ with respect to radiation source Sun  $S$, $\vec{V}$ the
corresponding velocity vector and $c$ the velocity of light. In the
expression of $F_{p}$, $L$ is luminosity of the radiating body,
while $m$, $\rho $ and $s$ are the mass, density and cross section
of the particle respectively.

The first term in equation (~\ref{eq:F}) expresses the radiation pressure. The second term represents the Doppler shift of the incident radiation and the third term is due to the absorption and subsequent re-emission of the incident radiation. These last two terms taken together are the Poynting-Robertson effect. The Poynting-Robertson effect will operate to sweep small particles of the solar system into the Sun at cosmically rapid rate.
 We consider the barycentric rotating co-ordinate system $Oxyz$ relative to inertial system with angular velocity $\omega$ and common $z$--axis.  We have taken line joining the primaries as $x$--axis. Let $m_1, m_2$ be the masses of bigger primary(Sun)  and smaller primary(Earth) respectively. Let  $Ox$, $Oy$  in the equatorial plane of smaller primary  and $Oz$ coinciding with the polar axis of $m_{2}$. Let $r_{e}$, $r_{p}$ be the equatorial and polar radii of $m_{2}$ respectively,  $r$ be the distance between primaries.  Let infinitesimal mass $m$ be placed at the  point $P(x,y,0)$. We take units such that sum of the masses and distance between primaries is  unity, the unit of time i.e. time period of $m_{1}$ about $m_{2}$  consists of $2\pi$ units such that the Gaussian constant of gravitational $\Bbbk^{2}=1$. Then perturbed mean motion $n$ of the primaries is given by $n^{2}=1+\frac{3A_{2}}{2}$, where $A_{2}=\frac{r^{2}_{e}-r^{2}_{p}}{5r^{2}}$ is oblateness coefficient of $m_{2}$.
Let $\mu=\frac{m_{2}}{m_{1}+m_{2}}$ then $1-\mu=\frac{m_{1}}{m_{1}+m_{2}}$ with $m_{1}>m_{2}$, where $\mu$ is mass parameter. Then coordinates of $m_{1}$ and $m_{2}$ are  $(-\mu,0)$ and $(1-\mu,0)$ respectively. Further, in our consideration, the velocity of light needs to be dimensionless, too, so consider  the dimensionless velocity of light as $c_d=c$ which depends on the physical masses of the two primaries and the distance between them. In this paper, we set  $c_d=299792458,\mu= 0.00003$ for all numerical results.

\begin{figure}
\plotone{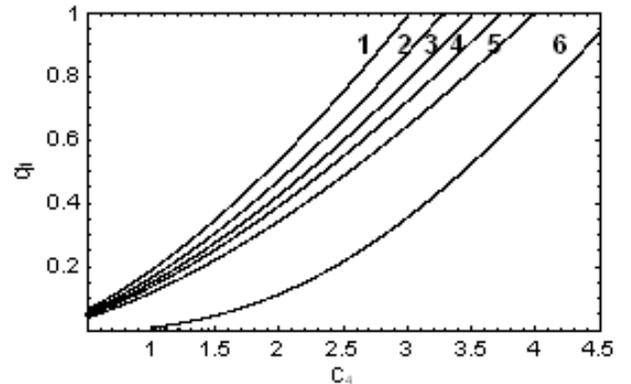}
\caption{The  Jacobi's constant $C_4-q_1$ for  $(1) A_2=0$ $(2) A_2=0.2 $ $(3) A_2=0.4$ $(4) A_2=0.6$ $(5) A_2=0.8$ $(6) A_2=1$, $0\leq q_1\leq 1$ \& $\mu= 0.00003$ }\label{fig:C4}
\end{figure}

\begin{figure}
\plotone{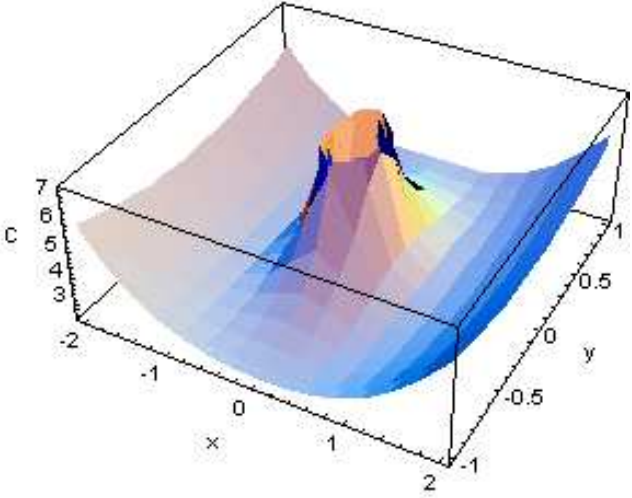}
\caption{$C$ Vs $x-y$ when  $A_2=0.0024, 0\leq q_1\leq 1$ \& $\mu= 0.00003$}
\label{fig:Cxy}
\end{figure}

\begin{figure}
\plotone{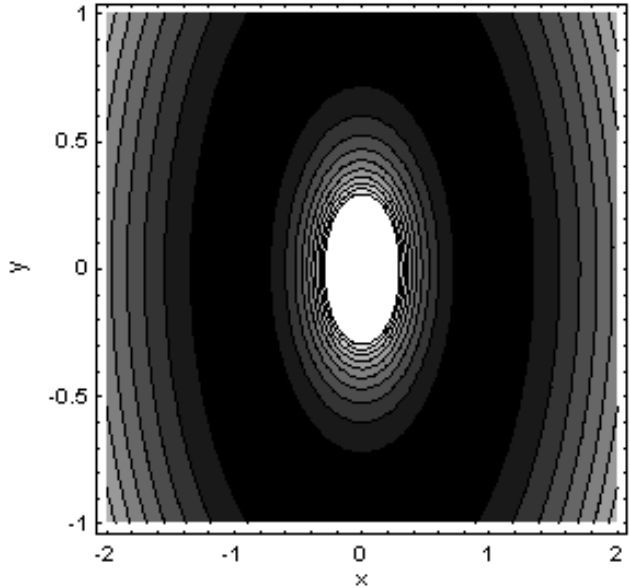}
 \caption{Contour plot shows  the  Jacobi's constant $C$ Vs $x-y$  when $A_2=0.0024, 0\leq q_1\leq 1$ \& $\mu=0.00003$}
\label{fig:CntCXY}
 \end{figure}

In the above mentioned reference system  the total   acceleration on the particle $P$ is as follows
\begin{eqnarray}
&&{\vec{a'}}=\vec{a}+2\vec{\omega}\times\vec{v}+\vec{\omega}\times(\vec{\omega}\times \vec{r})\nonumber\\&&=-\displaystyle{\frac{(1-\mu)\vec{r_1}}{r^3_1}-\frac{\mu{\vec{r_2}}}{r^3_2}-\frac{3}{2}\frac{\mu{A_2}\vec{r_2}}{r^5_2}}\nonumber\\
&&+\frac{(1-\mu)(1-q_1)}{r^2_1}\biggl\{\frac{\vec{r_1}}{r_1}-\frac{(\vec{\dot{r_1}}+\vec{\omega}\times\vec{r_1}).\vec{r_1}\vec{r_1}}{c_d r^2_1}\nonumber\\&&-\frac{\vec{\dot{r_1}}+\vec{\omega}\times{\vec{r_1}}}{c_d}\biggr\}\label{eq:E}
\end{eqnarray}
where
\begin{eqnarray*}
&&\vec{r_{1}}=x\hat{i}+y\hat{j},\
\vec{v}=\dot{x}\hat{i}+\dot{y}\hat{j},\
\vec{a}=\ddot{x}\hat{i}+\ddot{y}\hat{j},\
\vec{\omega}=n\hat{k},\\&&
\vec{r_{1}}=(x+\mu)\hat{i}+y\hat{j},
\vec{r_{2}}=(x+\mu-1)\hat{i}+y\hat{j},\\&&
r^{2}_{1}=(x+\mu)^{2}+y^{2},\
r^{2}_{2}=(x+\mu-1)^{2}+y^{2},\end{eqnarray*}
\begin{eqnarray*}
&&\vec\omega\times\vec v=-n(\dot y\hat i-\dot x \hat j),\,\vec \omega\times(\vec\omega\times\vec r)=n^2(x \hat i+y\hat j),\\&&
\vec{\dot{r_1}}+n\hat k\times\vec r_1=(\dot x -n y)\hat i +\left[ \dot y+n(x+\mu)\right]\hat j,\\&&
 \left(\vec{r_1}+n\hat k\times \vec r_1\right).\vec r_1 = \left[\dot x(x+\mu)+y\dot y \right]\end{eqnarray*}
 Substituting all these values in above relation (~\ref{eq:E}) we get as follows:

\begin{eqnarray*}
&& \left(\ddot x -2n \dot y-n^2x\right)\hat i +  \left(\ddot y +2n \dot x-n^2 y\right)\hat j
 \\&&=-\frac{(1-\mu)\vec r_1}{r_1^3}-\frac{\mu\vec r_2}{r_2^3}-\frac{3\mu\vec r_2A_2}{2r_2^5}+\frac{(1-\mu)q_1\vec r_1}{r_1^3}\\&&-\frac{(1-\mu)(1-q_1)}{r^2_1}\biggl\{\frac{(\vec{\dot{r_1}}+\vec{\omega}\times\vec{r_1}).\vec{r_1}\vec{r_1}}{c_d r^2_1}\\&&-\frac{\vec{\dot{r_1}}+\vec{\omega}\times{\vec{r_1}}}{c_d}\biggr\}\\&&=-\frac{(1-\mu)q_1\vec r_1}{r_1^3}-\frac{\mu\vec r_2}{r_2^3}-\frac{3\mu\vec r_2A_2}{2r_2^5}\\&&-\displaystyle{\frac{W_1}{c_d r^2_1}\biggl\{\frac{(\vec{\dot{r_1}}+\vec{\omega}\times\vec{r_1}).\vec{r_1}\vec{r_1}}{r^2_1}-\frac{\vec{\dot{r_1}}+\vec{\omega}\times{\vec{r_1}}}{1}\biggr\}}
\end{eqnarray*}
$W_1=\frac{(1-\mu)(1-q_1)}{c_d}$, substituting the values of $\vec r_1,\ \vec r_2$ and comparing the components of $\hat{i}$ and $\hat{j}$ we get the equations of motion of the infinitesimal mass particle in $x y$-plane.
\begin{eqnarray}
\ddot{x}-2n\dot{y}&=&U_{x} ,\label{eq:ux}\\
\ddot{y}+2n\dot{x}&=&U_{y} \label{eq:uy} \end{eqnarray}
where
\begin{eqnarray*}
&&U_x= n^{2}x-\frac{(1-\mu)q_1(x+\mu)}{r^3_1}-\frac{\mu(x+\mu-1)}{r^3_2}\\&&-\frac{3}{2}\frac{\mu{A_2}(x+\mu-1)}{r^5_2} \nonumber\\
&&-\frac{W_1}{r^2_1}\biggl\{\frac{(x+\mu)}{r^2_1}[(x+\mu){\dot{x}+y\dot{y}}] +\dot{x}-ny \biggr\},\\
&&U_y=n^{2}y
-\frac{(1-\mu)q_{1}{y}}{r^3_1}
-\frac{\mu{y}}{r^3_2}-\frac{3}{2}\frac{\mu{A_2}y}{r^5_2} \nonumber\\
&&-\frac{W_1}{r^2_1}\biggl\{\frac{y}{r^2_1}[(x+\mu)\dot{x}+y\dot{y}]+\dot{y}+n(x+\mu)\biggr\}\end{eqnarray*}
where
\begin{eqnarray}
&&U=\frac{n^2(x^2+y^2)}{2}+\frac{(1-\mu)q_1}{r_1}+\frac{\mu}{r_2}+\frac{\mu
 A_2}{2r_2^3}\nonumber\\&&+W_1\left\{ \frac{(x+\mu)\dot x + y\dot y}{2r_1^2}-n \arctan\left( \frac{y}{x+\mu}\right)\right\}\label{eq:FF}
 \end{eqnarray}
 $q_1=1-\frac{F_p}{F_g}$ is a mass reduction factor expressed in terms of the particle radius $\mathbf{a}$, density $\rho$ radiation pressure efficiency factor $\chi$ (in C.G.S. system): \( q_1=1-\frac{5.6\times{10^{-5}}}{\mathbf{a}\rho}\chi
\). The assumption $q_1=constant$ is equivalent to neglecting fluctuations in the beam of solar radiation and the effect of the planets shadow, obviously $q_1\leq1$. The energy integral of the problem is given by  $C=2U-{\dot{x}}^2-{\dot{y}}^2$, where the quantity $C$ is the Jacobi's  constant. The zero velocity curves[see (~\ref{fig:Cxy}),(~\ref{fig:CntCXY})]  are given by:
\begin{equation}
C_i=2U(x_i,y_i)\label{eq:C}
\end{equation}
Suffix $i=1,2,3,4,5$ correspond to respective $i^{th}$ Lagrangian equilibrium point $L_i$.
Using above relation we have determined $C_1\approx3.02978, C_2\approx 4.04133, C_3\approx3.53607$.  The values of $C_4(\approx C_5)$ are shown by different curves $(1-6)$ in figure (~\ref{fig:C4}) and table (~\ref{tbl-1}) for various values of $q_1, A_2$
\begin{table*}
\tabletypesize{\scriptsize}\caption{Jacobi's constant $C$ at  $L_4$}\label{tbl-1}
\begin{tabular}{crrrrr}
\tableline\tableline $A_2$ & $C_{4}:q_1=1$ & $C_{4}:q_1=0.75$& $C_{4}:q_1=0.5 $&$C_{4}:q_1=0.25$ & $C_{4}:q_1=0$\\
\tableline 
0.0&2.99997&2.47643&1.88988&1.19058&Indeterminate\\
0.2&3.27697&2.70632&2.06659&1.30303&Indeterminate\\
0.4&3.51014 &2.89885&2.21377&1.3961&Indeterminate\\
0.6&3.71975&3.0698&2.34214&1.47515&Indeterminate\\
0.8&3.98328&3.30328&2.53573&1.61041&Indeterminate\\
1&4.62593&4.06149&3.42503&2.62291&Indeterminate\\
\tableline
\end{tabular}
\end{table*}
\section{Existence of Equilibrium Points}
\label{sec:existEP}
What are they?
What are \lq\lq Lagrange points", also known as \lq\lq libration points" or \lq\lq L-points" or \lq\lq Equilibrium Points"?
These are all jargon for places where a light third body can sit  \lq\lq motionless" relative to two heavier bodies that are orbiting each other,  thanks to the force of gravity. The unstable Lagrange points --labeled $L_1, L_2$ and $L_3$ --lie along the line connecting the two large masses. The conditionally linearly stable Lagrange points in classical case are labeled $L_4$ and $L_5$ --form the apex of two equilateral triangles as  in figure (~\ref{fig:lagClass}), that have the large masses at their vertices.
\subsection{Collinear Equilibrium Points}
To investigate the Equilibrium Points divid the orbital plane $Oxy$ into three parts with respect to the primaries $x\leq-\mu$, $1-\mu\leq x$ and $-\mu<x<1-\mu$, $U_x=U_y=0$ then from equations (~\ref{eq:ux}) and (~\ref{eq:uy})
\begin{eqnarray*}
&&\biggl\{n^2-\frac{(1-\mu)q_1}{r^3_1}-\frac{\mu}{r^3_2}\bigl(1+\frac{3}{2}\frac{A_2}{r^2_2}\bigr)\biggr\}(x+\mu)+\frac{nW_1}{r^2_1}y\\&&=\mu\biggl\{n^2-\frac{1}{r^3_2}\bigl(1+\frac{3A_2}{2r^2_2}\bigr)\biggr\},\\
&&\biggl\{n^2-\frac{(1-\mu)q_1}{r^3_1}-\frac{\mu}{r^3_2}\bigl(1+\frac{3A_2}{2r^2_2}\bigr)\biggr\}y-\frac{nW_1}{r^2_1}(x+\mu)\\&&=0.
\end{eqnarray*}
\[\Rightarrow 
r_2^5\left(\frac{n W_1}{y\mu}-n^2\right)+r_2^2+\frac{3}{2}A_2=0,\]
\[
r_1^3\left[(1-\mu)yn^2+nW_1\right]-n(x+\mu)W_1r_1=(1-\mu)q_1y
\]
i.e.
\begin{eqnarray}
r_1=\left(\frac{q_1}{n^2}\right)^{1/3}\left[1-\frac{nW_1}{6(1-\mu)y}\right]\label{eq:r1},\\ 
r_2=\left[1-\frac{nW_1}{\mu y}(1-\frac{5}{2}A_2)]\right]^{-1/3} \label{eq:r2}
\end{eqnarray}
From above, we obtained:
\begin{eqnarray} &&x=-\mu\pm\nonumber\\&&\Bigl[\left(\frac{q_1}{n^2}\right)^{2/3}\left[1+\frac{nW_1}{2(1-\mu)y}+\frac{3A_2}{2}\right]^{-2/3}-y^2\Bigr]^{1/2} \label{eq:x1lag}\\
&&x=1-\mu\pm\Bigl[\left[1-\frac{nW_1}{\mu y}(1-\frac{5}{2}A_2)\right]^{-2/3}-y^2\Bigr]^{1/2}\label{eq:x2lag}
\end{eqnarray}
\begin{figure}
\plotone{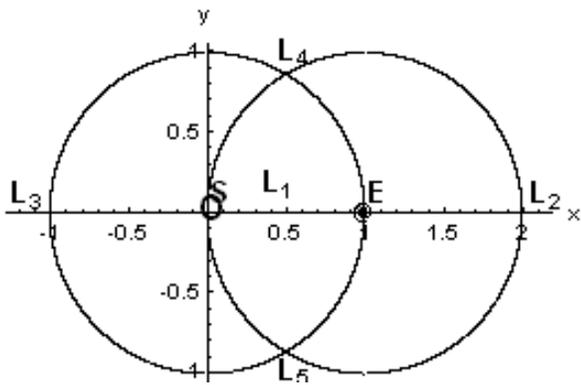}
\caption{Figure views the Lagrangian equilibrium points in classical case(when $A_2=0, q_1=1$,$\mu=0.00003$)}\label{fig:lagClass}
\end{figure}
\begin{figure}
\plotone{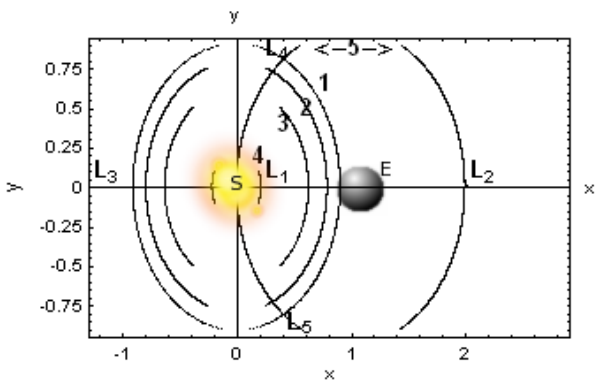}
\caption{The position of $L_1,L_2, L_3$, $L_4$ and $L_5$  for $(1) q_1=1$ $(2) q_1=0.75 $ $(3) q_1=0.25$ $(4) q_1=0$,  $(5) 0\leq q_1\leq 1$  when  $A_2=0.0024$ \& $\mu=0.00003$}\label{fig:lag12345}
\end{figure}
\begin{figure}
\plotone{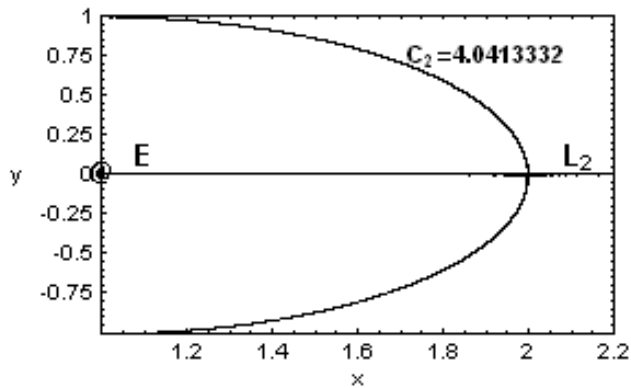}
\caption{Position of the  Lagrangian equilibrium point $L_2$  when  $A_2=0.0024$ \& $\mu=0.00003$, $0\leq q_1\leq 1$}\label{fig:lag2}
\end{figure}
\begin{figure}
\plotone{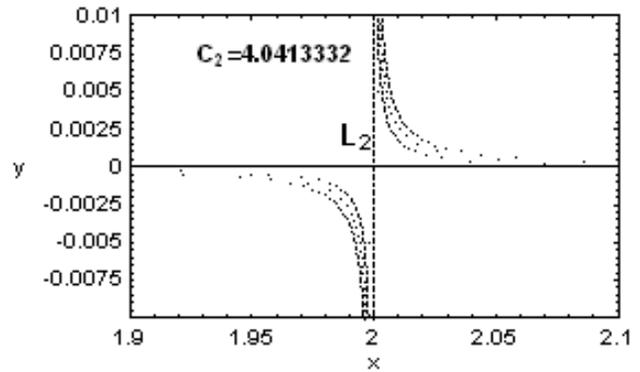}
\caption{Figure views  the  magnified region of $L_2$, when  $A_2=0.0024$ \& $\mu=0.00003$, $0\leq q_1\leq 1$}\label{fig:L2small}
\end{figure}
\begin{figure}
\plotone{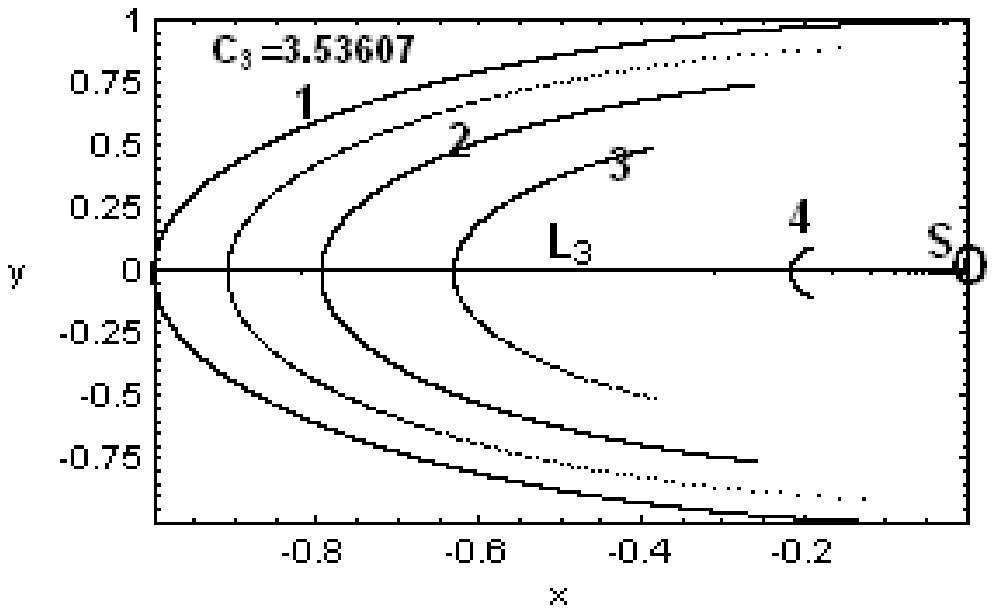}
\caption{Position of the  Lagrangian equilibrium points $L_3$  when  $A_2=0.0024$ \& $\mu=0.00003$ for  $(1) q_1=1$ $(2) q_1=0.75 $ $(3) q_1=0.25$ $(4) q_1=0$}\label{fig:lag3}
\end{figure}
\begin{figure}
\plotone{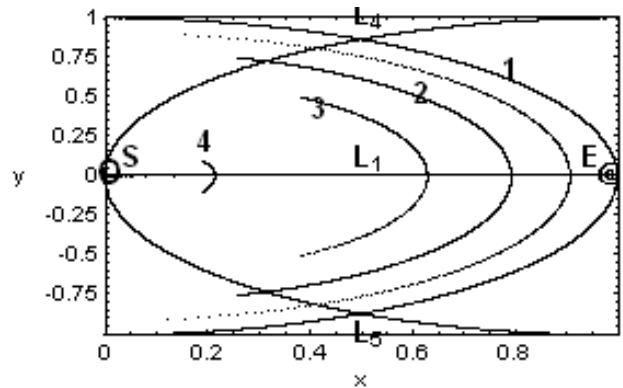}
\caption{Position of the  Lagrangian equilibrium points $L_1, L_4, L_5$ when  $A_2=0.0024$ \& $\mu=0.00003$, for  $(1) q_1=1$ $(2) q_1=0.75 $ $(3) q_1=0.25$ $(4) q_1=0$}\label{fig:lag145}
\end{figure}
Using equations (~\ref{eq:x1lag},~\ref{eq:x2lag}) the position of $L_1,L_2,L_3$ is presented graphically in  figures (~\ref{fig:lagClass},~\ref{fig:lag12345}) when $A_2=0.0024, 0\leq q_1\leq1$ \& $\mu=0.00003$. We have seen that the equilibrium points are no-longer collinear with the primaries. The values of $y$ are positive. We observe that the coordinates of $L_3$ are  the  functions of $q_1$  and $A_2$, corresponding   different curves labeled by $(1-4)$ are presented as in figures (~\ref{fig:lag3}). Now when $1-\mu\leq x$, there exists an equilibrium point $L_2$ this can be seen in figures  (~\ref{fig:lag2},~\ref{fig:L2small}). This  equilibrium point is also found away from the $Ox$ axis.  It is clear from the figure that $L_2$ is a function of $q_1$ and $A_2$. When $-\mu<x<1-\mu$, there are almost three equilibrium points they are $L_1,L_4$ and $L_5$ as in  figure  (~\ref{fig:lag145}).  It is clear from the figure, that the  $L_1$ is not collinear any more. $L_4$ is positioned above the $Ox$ axis while $L_5$ lies below it. All these results are similar with \citet*{Szebehely1967}, \citet*{Ragosetal1995}, \citet*{Papadakisetal2007Ap&SS} and others. 
\subsection{Triangular  Equilibrium Points}
\begin{table*}
\tabletypesize{\scriptsize}\caption{x co-ordinate of  $L_4$}\label{tbl-2}
\begin{tabular}{crrrrr}
\tableline\tableline $A_2$ & $x_{4}:q_1=1$ & $x_{4}:q_1=0.75$& $x_{4}:q_1=0.5 $&$x_{4}:q_1=0.25$ & $x_{4}:q_1=0$\\
\tableline 
0.0&{\bf 0.499949}&0.412692&0.314934&0.198383&Complex Infinity\\
0.25&0.374929&0.30949&0.236175&0.148765&Complex Infinity\\
0.5&0.249906&0.206284&0.157412&0.099146&Complex Infinity\\
0.75&0.12488&0.103076&0.0786479&0.0495247&Complex Infinity\\
1&-0.000148407&-0.000134579&-0.000118667&-0.0000980564&Complex Infinity\\
\tableline
\end{tabular}
\end{table*}
\begin{table*}
\tabletypesize{\scriptsize}\caption{y co-ordinate of  $L_4$}\label{tbl-3}
\begin{tabular}{crrrrr}
\tableline\tableline $A_2$ & $y_{4}:q_1=1$ & $y_{4}:q_1=0.75$& $y_{4}:q_1=0.5 $&$y_{4}:q_1=0.25$ & $y_{4}:q_1=0$\\
\tableline 
0.0&{\bf 0.86605}&0.809425&0.728552&0.597927&Indeterminate\\
0.25&0.847795&0.790465&0.709777&0.581035&Indeterminate\\
0.5&.790543&0.73068&0.650249&0.527174&Indeterminate\\
0.75&0.684579&0.61834&0.536576&0.422435&Indeterminate\\
1&0.499831&0.412542&0.314727&0.198036&Indeterminate\\
\tableline
\end{tabular}
\end{table*}
For the triangular equilibrium points $y\neq{0}$, $U_x=U_y=0$, then from equations (~\ref{eq:ux}) and (~\ref{eq:uy}) we get as follows:
\begin{figure}
\plotone{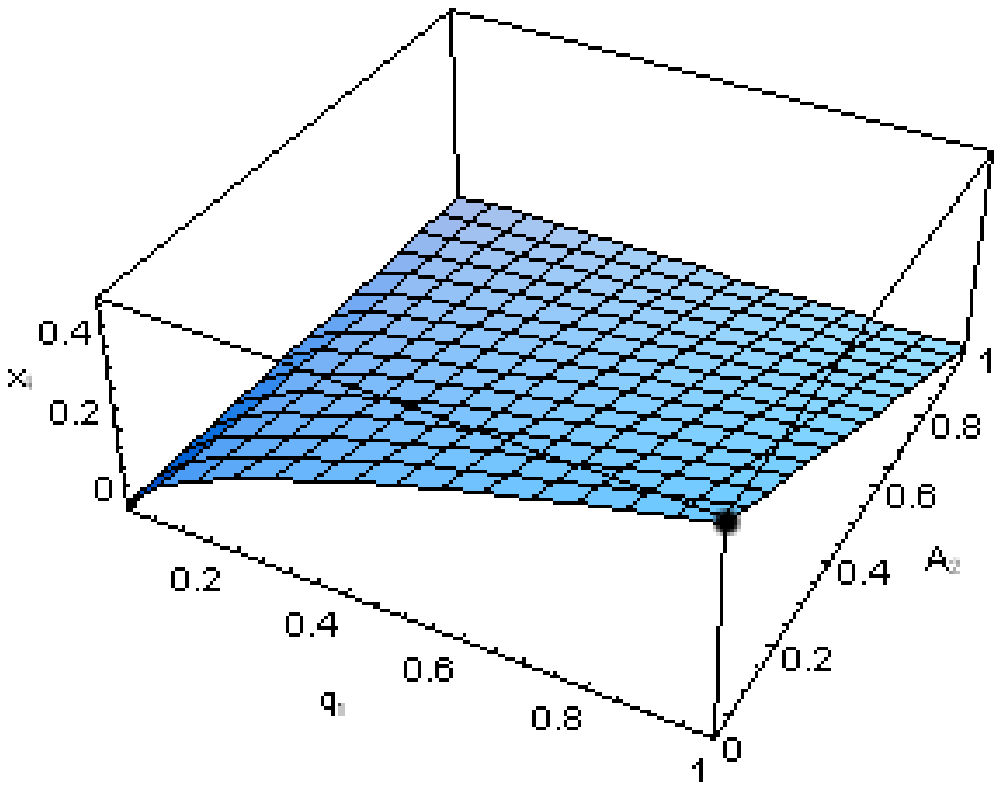}
\caption{The  $x_4$ coordinate of equilibrium point when $0\leq A_2\leq1, 0\leq q_1\leq 1$ \& $\mu= 0.00003$, black dot on figure  indicates the  value of $x_4$ in classical case} \label{fig:x}
 \plotone{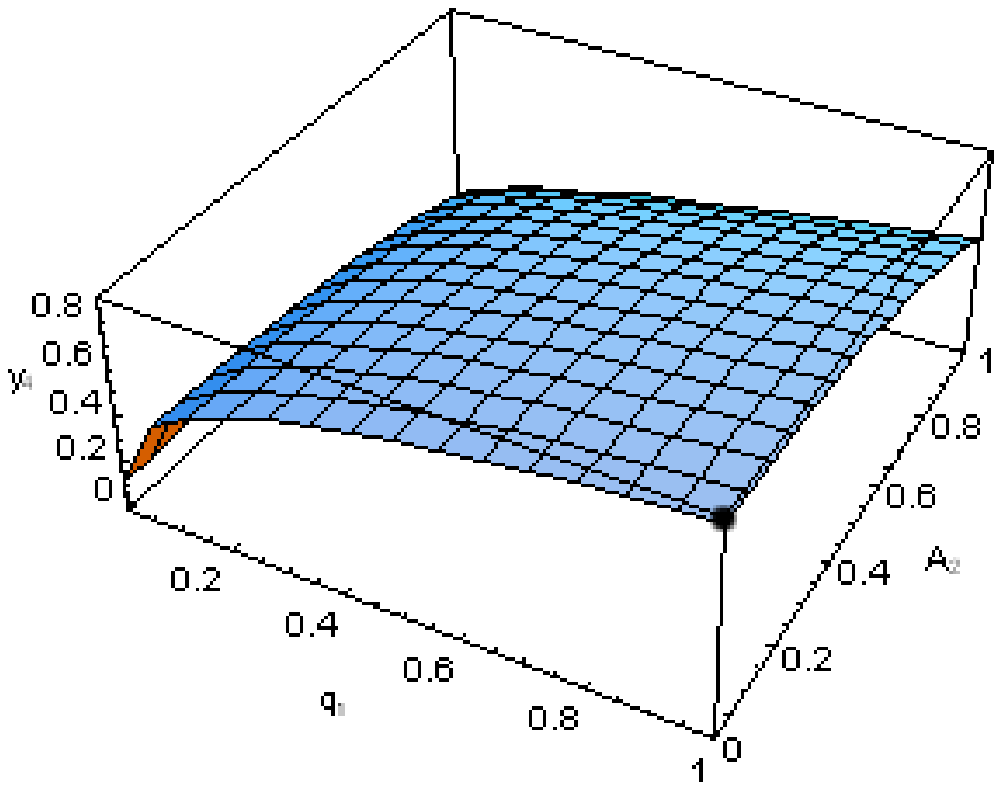}
\caption{The  $y_4$ coordinate of equilibrium point when  $0\leq A_2\leq1, 0\leq q_1\leq 1$ \& $\mu= 0.00003$, black dot on figure  indicates the  value of $y_4$ in classical case}\label{fig:y}
\end{figure}
\begin{eqnarray}
&&n^{2}x-\frac{(1-\mu){q_1}(x+\mu)}{r^3_1}-\frac{\mu(x+\mu-1)}{r^3_2}\nonumber\\&&-\frac{3}{2}\frac{\mu{A_2}(x+\mu-1)}{r^5_2}+\frac{W_1ny}{r^2_1}=0,\label{eq:ux0}\\&&
n^{2}y-\frac{(1-\mu){q_1}y}{r^3_1}-\frac{\mu{y}}{r^3_2}-\frac{3}{2}\frac{\mu{A_2}y}{r^5_2}\nonumber\\&&-\frac{W_1}{r_1^2}n(x+\mu)=0.\label{eq:uy0}
\end{eqnarray}
\begin{eqnarray*}
&&\biggl\{n^2-\frac{(1-\mu)q_1}{r^3_1}-\frac{\mu}{r^3_2}\bigl(1+\frac{3}{2}\frac{A_2}{r^2_2}\bigr)\biggr\}(x+\mu)+\frac{nW_1}{r^2_1}y\nonumber\\&&=\mu\biggl\{n^2-\frac{1}{r^3_2}\bigl(1+\frac{3A_2}{2r^2_2}\bigr)\biggr\}, \label{eq:x1}\nonumber\\&&
\biggl\{n^2-\displaystyle{\frac{(1-\mu)q_1}{r^3_1}-\frac{\mu}{r^3_2}\bigl(1+\frac{3A_2}{2r^2_2}\bigr)}\biggr\}y-\frac{nW_1}{r^2_1}(x+\mu)\nonumber\\&&=0.
\end{eqnarray*}
From  equations (~\ref{eq:ux0}) and  (~\ref{eq:uy0}),  we obtained as follows:
\begin{equation}
\biggl\{n^2-\displaystyle{\frac{1}{r^3_2}\bigl(1+\frac{3A_2}{2r^2_2}\bigr)}\biggr\}\mu{y}=nW_1.
\end{equation}
In the case of photogravitational restricted three body problem we have $r_{1_0}=q_1^{1/3}=\delta$ and $r_{2_0}=1$. We suppose that due to P-R drag and oblateness, perturbations in $r_{1_0}$, $r_{2_0}$ are $\epsilon_1$, $\epsilon_2$ respectively,  where ${\epsilon_{i}}'$s are very small quantities, then
\begin{equation}
r_1=q^{1/3}_1(1+\epsilon_1),\quad r_2=1+\epsilon_2
\end{equation}
Putting these values in equation (~\ref{eq:uy0}) and neglecting higher order terms of small quantities, we get
\begin{eqnarray*}
\left[ n^2-1+3\epsilon_2-\frac{3}{2}A_2+\frac{15}{2}A_2 \epsilon_2\right]\mu y=nW_1,\\
3\epsilon_2\left(1+\frac{5}{2}A_2\right)\mu y=nW_1,
\end{eqnarray*}
\[
\epsilon_2=\displaystyle{\frac{nW_1(1-\frac{5}{2}A_2)}{3\mu{y_0}}},
\]
where $y_0$ is y-coordinate of triangular equilibrium points in  photogravitational restricted three body problem case. Using  values  $r_1$,  $r_2$ with equation (~\ref{eq:uy0}), we get  following:
\begin{eqnarray*}
&&\biggl\{n^2-\frac{(1-\mu)q_1}{r^3_1}-\frac{\mu}{r^3_2}\bigl(1+\frac{3A_2}{2r^2_2}\bigr)\biggr\}y-\frac{nW_1}{r^2_1}(x+\mu)\\&&=0,\\
&&\biggl[\biggl\{n^2-(1-\mu)(1-\epsilon_1)\\&&-\mu(1-\epsilon_2)\bigl(1+\frac{3A_2(1-2\epsilon_2)}{2}\bigr)\biggr\}\biggr]y\\&&-\frac{nW_1(r_1^2-r_2^2+1)}{r^2_1}=0,\\&&
\left[ n^2-(1-\mu)+3(1-\mu)\epsilon_1-\mu\right.\nonumber\\&&+\left.3\epsilon_2-\mu(1-5\epsilon_2)\frac{3}{2}A_2\right] y=\frac{nW_1}{2}
\end{eqnarray*}
we obtained \(\epsilon_1=-\displaystyle{\frac{nW_1}{6(1-\mu){y_0}}}-\frac{A_2}{2}
\),  \(r_1=\delta\biggl\{1-\frac{nW_1}{6(1-\mu){y_0}}-\frac{A_2}{2}\biggr\}\), \( r_2=1+\frac{nW_1}{3\mu{y_0}}(1-\frac{5}{2}A_2)
 \).

Now we have  $x+\mu=\frac{r^2_1-r^2_2+1}{2}$,\ $y^2=r^2_1-(x+\mu)^2$ therefore,
\begin{eqnarray}
&&(x+\mu)=\frac{1}{2} \left[\delta^2\biggl\{1- \frac{nW_1}{6\amc {y_0}}-\frac{A_2}{2}  \biggr\}^2\right.\nonumber\\&&\left.-\biggl\{1+\frac{nW_1}{3\mu{y_0}\abc}\biggr\}^2+1\right]\nonumber
\end{eqnarray}
or
\begin{eqnarray}
&&x_4=x_0\biggl[1-\adc\frac{A_2}{x_0}\nonumber\\&&-\displaystyle{\frac{nW_1[\amc\abc+\mu\aac\adc]}{3\mu\amc y_0 x_0}}\biggr]\label{eq:xl4}\end{eqnarray}
\begin{eqnarray}
&&y^2=r^2_1-(x+\mu)^2 \nonumber \\
 &=&\delta^2\biggl[1-\frac{nW_1}{6\amc y_0}\biggr]^2-\biggl[\adc(1-\adc A_2)\nonumber\\&&-\frac{nW_1[\amc\abc+\mu\aac\adc]}{3\mu\amc y_0}\biggr]^2 \nonumber 
 \end{eqnarray}
or
\begin{eqnarray}
 &&y_4=y_0 \biggl[1-\frac{\delta^2(1-\adc)A_2}{y^2_0}\nonumber\\&&-\frac{nW_1\delta^2[2\mu-1-\mu(1-\frac{3A_2}{2})\adc+7\amc\frac{A_2}{2}]}{3\mu\amc y^3_0}\biggr]^{1/2}\label{eq:yl4} 
 \end{eqnarray}
where $(x_0,y_0)$ are coordinates of $L_4,L_5$ in  photogravitational restricted three body problem case as:
\[
x_0=\adc-\mu,\quad y_0=\pm\delta\biggl(1-\frac{\delta^2}{4}\biggr)^{1/2}, \quad \delta=q^{1/3}_1
\]
The position of $L_{4(5)}$ is given by equations (~\ref{eq:xl4}), (~\ref{eq:yl4}) which are valid for $W_1\ll 1,A_2\ll 1$. For simplicity we suppose $\gamma=1-2\mu$,  $q_1=1-\epsilon$, with $|\epsilon|<<1$ then the coordinates $(x_4,y_4)$ of $L_{4(5)}$  can be written as follows:
\begin{eqnarray}
x_4&=&\frac{\gamma}{2}-\frac{\epsilon}{3}-\frac{A_2}{2}+\frac{A_2
\epsilon}{3}\nonumber\\&&-\frac{(9+\gamma)}{6\sqrt{3}}nW_1-\frac{4\gamma
\epsilon}{27\sqrt{3}}nW_1 \\
y_4&=&\frac{\sqrt{3}}{2}\Bigl\{1-\frac{2\epsilon}{9}-\frac{A_2}{3}-\frac{2A_2
\epsilon}{9}\nonumber\\&&+\frac{(1+\gamma)}{9\sqrt{3}}nW_1-\frac{4\gamma
\epsilon}{27\sqrt{3}}nW_1\Bigr\}
\end{eqnarray}

Figures (~\ref{fig:x},~\ref{fig:y}) and  tables (~\ref{tbl-2},~\ref{tbl-3}) show how the coordinates $(x_4,y_4)$  points are decreasing functions of $A_2, q_1, W_1$. Black dots on the figures indicate the $(x_4,y_4)=\left(\frac{1}{2}-\mu, \frac{\sqrt{3}}{2}\right)$.

\section{Comments on the Linear  Stability}
\label{sec:cmnt_lstb}
In order to study the linear stability of any Lagrangian equilibrium point $L_i(i=1-5)$ the origin of the coordinate system to its position $(x_i,y_i)$ by means of  $x=x_i+\alpha$,\,  $y=y_i+\beta$, where  $\alpha=\xi e^{\lambda{t}}$,\ $\beta=\eta e^{\lambda{t}}$ are the small displacements  $\xi,\eta$,\  $\lambda$ these parameters, have to be determined. 
 Therefore the equations of perturbed motion corresponding to the system of equations (~\ref{eq:ux}), (~\ref{eq:uy}) may be written as follows:
\begin{align}
\ddot{\alpha}-2n\dot{\beta} &= {\alpha}{U^i_{xx}}+{\beta}{U^i_{xy}}+\dot{\alpha}{U^i_{x\dot{x}}}+{\dot{\beta}}{U^i_{x\dot{y}}} \\
\ddot{\beta}+2n\dot{\alpha}&= {\alpha}{U^i_{yx}}+{\beta}{U^i_{yy}}+\dot{\alpha}{U^i_{y\dot{x}}}+\dot{\beta}{U^i_{y\dot{y}}}
\end{align}
where superfix $i$ is corresponding to the $L_i(i=1-5)$
 \begin{align}(\lambda^2-\lambda{U^i_{x\dot{x}}}-{U^i_{xx}})\xi
+[-(2n+{U^i_{x\dot{y}}})\lambda-{U^i_{xy}}]\eta&=0\label{eq:lambda_x}\\
[(2n-{U^i_{y\dot{x}}})\lambda-{U^i_{yx}}]\xi
+(\lambda^2-\lambda{U^i_{y\dot{y}}}-{U^i_{yy}})\eta&=0\label{eq:lambda_y}
\end{align}
Now  above system has singular solution if,
\[
\begin{vmatrix}
\lambda^2-\lambda{U^i_{x\dot{x}}}-U^i_{xx}& -(2n+{U^i_{x\dot{y}}})\lambda-U^i_{xy} \\(2n-{U^i_{y\dot{x}}})\lambda-U^i_{yx}& \lambda^2-\lambda{U^i_{y\dot{y}}}-U^i_{yy}\\
\end{vmatrix}
=0
\]
\begin{equation}
\Rightarrow \quad \lambda^4+a\lambda^3+b\lambda^2+c\lambda+d=0 \label{eq:cheq}
\end{equation}
At the equilibrium points equations (~\ref{eq:ux}),(~\ref{eq:uy}) gives us the following:
\begin{eqnarray*}
&&a=3\zws,\, b=2n^2-f_i-\frac{3\mu{A_2}}{{r^5_2}_i}+\frac{2W^2_1}{{r^4_1}_i}\\
&&c=-a(1+e),\\&& e=\frac{\mu}{{r^5_2}_i}{A_2}+\frac{\mu}{{{r^2_1}_i}{{r^5_2}_i}}\zabs{y^2_i},\\&&
 d=(n^2-f_i)\biggl[n^2+2f_i-\frac{3\mu{A_2}}{{r^5_2}_i}\biggr]\nonumber\\&&+\frac{9\mu\amc{q_1}}{{{r^5_1}_i}{{r^5_2}_i}}\zabs{y^2_i}\nonumber\\&&
-\zd{\frac{6nW_1}{{{r^4_1}_i}{{r^5_2}_i}}\zabs\Bigl\{(x_i+\mu)(x_i+\mu-1)+y^2_i\Bigr\}},\\&&f_i=\frac{(1-\mu)q_1}{r^3_{1_i}}+\frac{\mu}{r^3_{2i}}\bigl(1+\frac{3}{2}\frac{A_2}{r^2_{2i}}\bigr) \end{eqnarray*}

The points $L_1,L_2,L_3$  no longer lie along the line joining the primaries, since the condition is not satisfied for them, so taking $y\rightarrow 0,\frac{W_1}{y}\rightarrow 0$ because $y>>W_1, x>>W_1$, from (~\ref{eq:r1}) we have  $r_1\approx \frac{q_1^{1/3}}{n^2}$. In this case  $f_i>1$ for $i=1,2,3$, so that for each  $L_1,L_2,L_3$  characteristic equation (~\ref{eq:cheq}) has at least one positive root, this implies that these points are unstable in the sense of Lyapunov.
\subsection{Linear Stability without P-R Effect}
\label{subsec:lstb_No_PR}
 Now we consider the problem when P-R effect is not included($W_1=0$), then $r_1=q_1^{1/3}(1-\frac{A_2}{2}), r_2=1$. The coordinates of triangular points $L_{4(5)}$ are as:
\begin{eqnarray}
x&=&\left(\frac{q_1^{1/3}}{2}-\mu\right)-\frac{q_1^{1/3}A_2}{2}\label{eq:xnpr}\\
y&=&q_1^{1/3}\biggl[\left(1-\frac{q_1^{2/3}}{4}\right)^{1/2}\nonumber\\&&-\left(1-\frac{q_1^{2/3}}{2}\right)\left(1+\frac{q_1^{2/3}}{4}\right)A_2\biggr]\label{eq:ynpr}
\end{eqnarray}
\begin{eqnarray*}
&&a=0, c=0,f_i=n^2,(i=4,5)\\
&&b=n^2-3\mu{A_2}, d=9\mu\amc g
\end{eqnarray*}
where \(g=(1-A_2)\biggl[{1-\frac{q_1^{2/3}(1-A_2)}{4}}\biggr]\). From characteristic equation (~\ref{eq:cheq}) we obtained,
\begin{equation}
 \lambda^2=\frac{-b\pm{(b^2-4d)}^{1/2}}{2}\label{eq:cheqNpr}
\end{equation}
For stable motion  $0<4d<b^2$, i.e.
\begin{equation*}
{\bigl(n^2-3\mu{A_2}\bigr)^2}>36\mu\amc g
\end{equation*}
In classical case $A_2=0$, $q_1=1$, $W_1=0$, $n=1$, we have following:
\(
1>27\mu\amc\  \Rightarrow \quad  \mu <0.0385201.
\)
The possible roots of equation (~\ref{eq:cheqNpr}) are given in table(~\ref{tbl-6}), we see that all the roots are purely imaginary quantities. Hence the triangular  equilibrium points are stable in the sense of Lyapunov stability provided  $\mu< \mu_{Routh}=0.0385201$.

Equation (~\ref{eq:cheqNpr}) has imaginary roots $\lambda_{1,2}=\pm \mathbf{i}\omega_1$, $\lambda_{3,4}=\pm \mathbf{i}\omega_2$, $\mathbf{i}=\sqrt{-1}$,  this gives us:
\begin{equation}\omega_{1,2}=\left(\frac{-b\pm{(b^2-4d)}^{1/2}}{2}\right)^{1/2}\label{eq:omega12}\end{equation}

There are  three main  cases  of resonances as:
\begin{equation}
\omega_1-k\omega_2=0,\quad k=1,2,3\label{eq:Reson}
\end{equation}
For $k=1$ we have positive stable resonance and for $k=2,3$ we have unstable resonances. Using (~\ref{eq:omega12}) and (~\ref{eq:Reson}) we obtained a root of mass parameter.
\begin{eqnarray}
\mu_k=\frac{3g+2 K A_2+3K A_2^2-\sqrt{g}\sqrt{9g-4K+9K A_2^2}}{6(g+K A_2^2)}
\end{eqnarray}
where $K=\frac{k^2}{(k^2+1)^2}$.
Now we suppose  $q_1=1-\epsilon$, with $|\epsilon|<<1$, neglecting higher  order terms, we obtained the critical mass parameter values corresponding to $k=1,2,3$ as :
\begin{eqnarray}
\mu_1&=&0.0385208965 + 
    0.6755841373 A_2\nonumber\\&&- 
    0.0089174706\epsilon\\
\mu_2&=& 0.0242938971+ 
    0.4322031625A_2\nonumber\\&&- 
    0.0055364958\epsilon\\
\mu_3&=&0.0135160160+ 
    0.2430452832A_2\nonumber\\&&- 
    0.0030452832\epsilon
\end{eqnarray}
\begin{figure}[h]
\plotone{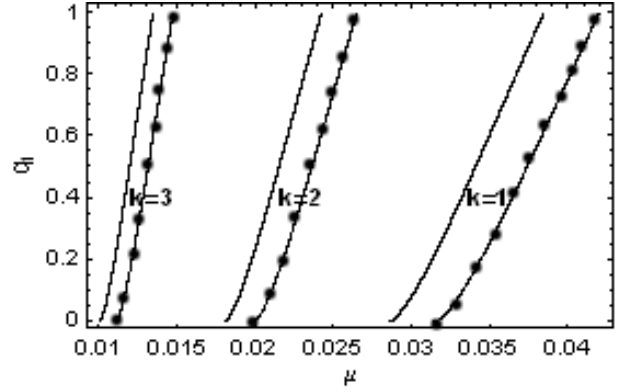}
\caption{The linear stability region  for $\mu-q_1$ parameter space and resonance curves $\omega_1-k\omega_2=0, k=1,2,3$. at $A_2=0, 0.02$(doted lines) \& $\mu=0.00003$,  $0\leq q_1\leq 1$}\label{fig:muk}
\end{figure}
  The linear stability region and corresponding main resonance curves in $\mu-q_1$ parameter space are shown in figure (~\ref{fig:muk}) the doted lines are  corresponding  to $A_2=0.02$, the curve corresponding to $k=1,(q_1=1, A_2=0,\mu_1=\mu_{Routh}=0.038521)$ is actual boundary of the stability region these results are similar to  \citet*{Markellos1996Ap&SS} and others.
The critical values of mass parameter $\mu$ are given in the tables (~\ref{tbl-4},~\ref{tbl-5}) for various values of $q_1, A_2$. The classical critical values of $\mu$ are similar to \citet*{Deprit1967}. We observe that the effect of radiation pressure reduces the  linear stability zones,  these are also affected by oblateness of second primary.
\begin{table*}
\tabletypesize{\scriptsize}\caption{$A_2=0$}\label{tbl-4}
\begin{tabular}{crrrrr}
\tableline\tableline $k$ & $\mu_{k}:q_1=1$ & $\mu_{k}:q_1=0.75$  & $\mu_{k}:q_1=0.5 $&$\mu_{k}:q_1=0.25$ &$\mu_{k}:q_1=0$ \\
\tableline 
1&0.0385209 & 0.0363201& 0.0341355 & 0.0318518&0.0285955\\
2&0.0242939 &0.0229262& 0.0215661 & 0.0201415&0.0181056\\
3& 0.013516 & 0.0127632& 0.0120136&0.0112275& 0.0101021\\
4&0.00827037&0.0078121&0.00735548&0.00687629&0.00618979\\
5&0.0055092 &0.00520474&0.004901287&0.00490128&0.00412616\\
\tableline
\end{tabular}
\end{table*}
\begin{table*}
\tabletypesize{\scriptsize}\caption{$ A_2=0.02$}\label{tbl-5}
\begin{tabular}{crrrrr}
\tableline\tableline $k$ & $\mu_{k}:q_1=1$ & $\mu_{k}:q_1=0.75$& $\mu_{k}:q_1=0.5 $&$\mu_{k}:q_1=0.25$ & $\mu_{k}:q_1=0$\\
\tableline 
1&0.0413469 &0.0390477&0.0367581&0.0343567&0.0309186\\
2&0.026094&0.0246633&0.0232361&0.0217366&0.019585\\
3&0.0145252 &0.0137369&0.0129496&0.0121214&0.0109312\\
4&0.00889015&0.00841007&0.00793029&0.00742525&0.00669897\\
5&0.00592287&0.00560384&0.00528491&0.00494909&0.00446598\\
\tableline
\end{tabular}

\end{table*}
\subsection{Linear Stability with P-R Effect}
\label{subsec:lstb_with_PR}
Now consider the  problem when P-R Effect is included i.e., $q_1\neq1,W_1\neq 0$, $A_2\neq0$. Using Ferrari's theorem the roots of characteristic equation (~\ref{eq:cheq}) are given by:
\begin{equation}
\lambda_i=-\frac{(a+A)}{4}\pm\sqrt{\left(\frac{a+A}{4}\right)^2-B}\label{eq:cheq_withPR}
\end{equation}
where $A=\pm\sqrt{8l-4b+a^2}$, $B=l(1+\frac{a}{A})-\frac{c}{A}$, $ i=1,2,3,4$ and  $l$ is any real root of the equation
\begin{eqnarray*}
8l^3-4bl^2+(2ac-8d)l+d(4b-a^2)-c^2=0\label{eq:l}
\end{eqnarray*} 
This can be written as:
\begin{equation}
2l^2-bl^2-2dl+db=\frac{a^2}{4}\left\{(1+e)^2-2(1+e)l+d\right\}\label{eq:Orl}
\end{equation}
This equation has an exact real root $l=\frac{b}{2}$ for $a=0$. When $a\neq0$ the roots of characteristic equation (~\ref{eq:cheq}) will be obtained  in form of the rapidly convergent series
\begin{equation}
l=\frac{b}{2}+\sum_{j=1}^\infty \alpha_j a^{2j}\label{eq:ser_l}
\end{equation}
Using (~\ref{eq:Orl}, ~\ref{eq:ser_l}), taking the coefficients of $a^2$ only, we get $A=a\pm\sqrt{1+8\alpha_1}$, 
\[ B=\left( \frac{b}{2}+\alpha_1a^2\right)\left(1\pm\sqrt{1+8\alpha_1}\right)\mp\frac{1+e}{\sqrt{1+8\alpha_1}},\]
\begin{equation}
\alpha_1=\frac{(1+e)(1+e^2-b)+d}{2(b^2-4d)}>0
\end{equation}
We obtained the characteristic roots
\begin{eqnarray}
\lambda_{1,2}&=&-\frac{a\left(1+\sqrt{1+8\alpha_1}\right)}{4}\nonumber\\&&\pm\sqrt{\frac{a^2\left(1+\sqrt{1+8\alpha_1}\right)}{16}-B_1}\label{eq:lambda12}\\
\lambda_{3,4}&=&-\frac{a\left(1-\sqrt{1+8\alpha_1}\right)}{4}\nonumber\\&&\pm\sqrt{\frac{a^2\left(1-\sqrt{1+8\alpha_1}\right)}{16}-B_2}\label{eq:lambda34}
\end{eqnarray}
where
\[ B_1=\left( \frac{b}{2}+\alpha_1a^2\right)\left(1+\sqrt{1+8\alpha_1}\right)-\frac{1+e}{\sqrt{1+8\alpha_1}},\]
\[ B_2=\left( \frac{b}{2}+\alpha_1a^2\right)\left(1-\sqrt{1+8\alpha_1}\right)+\frac{1+e}{\sqrt{1+8\alpha_1}}\]
Using above equations(~\ref{eq:lambda12},~\ref{eq:lambda34}), we obtained the roots of characteristic  equation (~\ref{eq:cheq}) which  are presented in table (~\ref{tbl-7}) for various values of $q_1, A_2$. We see that at least one of the roots $\lambda_i(i=1,2,3,4)$  have a positive real part due to P-R effect  as in \citet*{Chernikov1970}. Hence the triangular  equilibrium points are unstable in the sense of Lyapunov stability.
\section{Conclusion}
\label{sec:con}
The distances of $L_1, L_2$ from the second primary are monotonically increasing with $\mu$ . The Jacobi's constant at the $L_1, L_3$ increases monotonically with $\mu$.  If $C<C_1$, the particles can leave the system while if  $C_1<C<C_2$, i.e. the third body is not confined to its motion around the Sun but it is allowed to become a satellite of the Earth. This does not mean that it will become permanently a satellite of the Earth since its Jacobi's constant is such that it is not confined to the zero velocity oval around Earth. The particle with $C>C_2$ then it can not change position from the vicinity of Sun to that of Earth. The particles with very high $C$ values have low relative energy levels and they either move around one of the primaries or move for outside of the system. Finally we have conclude that, if one the primary is exerts the light radiation pressure and second primary is an oblate spheroid, then the  gravitational radiation force and oblateness influence the existence, location of equilibrium points. In classical case when $q_1$, $A_2=0$, we have the collinear points $(L_1,L_2,L_3)$ and two triangular equilibrium points $L_4,L_5$.  We have seen that $L_1,L_2,L_3$ are linearly unstable, while  $L_4,L_5$ are conditionally stable in the sense of Lyapunov when P-R effect is not considered. The effect of radiation pressure reduces the  linear stability zones while P-R effect induces an instability in the sense of Lyapunov.

\acknowledgments{I am very thenkful to Dr. Uday Dolas, Dr. Bhola Ishwar and  Mr. V.S. Srivastava for their persuasion. I am also  thankful for the referees\rq\  comments and suggestions and they have been very useful in improving the manuscript.}
\clearpage
\begin{deluxetable}{rrrrrrrrr}
\tabletypesize{\scriptsize}
\rotate
\tablecaption{Roots of characteristic equation when P-R-effect is not considered \label{tbl-6}}
\tablewidth{0pt}
\tablehead{
\colhead{$q_1$} & \colhead{$\lambda_{1,2}:A_2=0$} & 			\colhead{$\lambda_{3,4}:A_2=0$} & \colhead{$\lambda_{1,2}:A_2=.02$} & \colhead{$\lambda_{3,4}:A_2=.02$} &
\colhead{$\lambda_{1,2}:A_2=.04$} & \colhead{$\lambda_{3,4}:A_2=.04$}&
\colhead{$\lambda_{1,2}:A_2=.06$} & \colhead{$\lambda_{3,4}:A_2=.06$}
}
\startdata
1 &$ \pm 0.0100624\mathbf{i}$ &$\pm1.414178\mathbf{i}$&$\pm0.0103471\mathbf{i} $&$\pm0.0103471\mathbf{i}$&$\pm1.435232\mathbf{i}$&$\pm1.45598\mathbf{i}$&$\pm0.0109142\mathbf{i}$&$\pm1.476440\mathbf{i}$\\

0.8 & $\pm0.0102916 \mathbf{i}$ & $\pm1.414176\mathbf{i}$&$\pm0.010577 \mathbf{i}$& $\pm1.435230\mathbf{i}$&$\pm 0.0108624\mathbf{i}$& $\pm1.45598\mathbf{i}$&$\pm0.0111478\mathbf{i}$&$\pm1.476438\mathbf{i}$\\

0.6& $\pm0.0105354\mathbf{i} $&$\pm1.414174\mathbf{i}$ &$\pm0.0108217$&$\pm1.435229\mathbf{i}$&$\pm0.0111086\mathbf{i}$&$\pm1.455978\mathbf{i}$&$\pm0.0113964\mathbf{i}$&$\pm1.476436\mathbf{i}$\\

0.4&$\pm0.0108019 \mathbf{i}$&$\pm1.414172\mathbf{i}$&$\pm 0.0110893\mathbf{i} $&$\pm1.435226\mathbf{i}$&$\pm.3704\mathbf{i}$&$\pm1.455976\mathbf{i}$&$\pm0.0116684 \mathbf{i}$&$\pm1.476434\mathbf{i}$\\

0.2&$\pm0.0347\mathbf{i}$&$\pm1.414170\mathbf{i}$&$\pm0.0114002\mathbf{i}$&$\pm1.435224\mathbf{i}$&$\pm0.011378\mathbf{i}$&$\pm  1.455973\mathbf{i}$&$\pm0.0119841\mathbf{i}$&$\pm1.476432\mathbf{i}$\\

0.0 &&$\mathbf{i}$&$1.435269\mathbf{i}$&0&$1.456020\mathbf{i}$&0& $\pm1.476480\mathbf{i}$&0\\
\enddata

\tablecomments{Table \ref{tbl-6}  presents the roots of characteristic equation(~\ref{eq:cheq}) for $\mu=0.00003$, when P-R effect is not considered}
\end{deluxetable}

\begin{deluxetable}{rrrrrrr}
\tabletypesize{\scriptsize}
\rotate
\tablecaption{Roots of characteristic equation when P-R effect is considered\label{tbl-7}}
\tablewidth{0pt}
\tablehead{
\colhead{$q_1$} & \colhead{$\lambda_{1,2}:A_2=0$} & 			\colhead{$\lambda_{3,4}:A_2=0$} & \colhead{$\lambda_{1,2}:A_2=.02$} & \colhead{$\lambda_{3,4}:A_2=.02$} &
\colhead{$\lambda_{1,2}:A_2=.04$} & \colhead{$\lambda_{3,4}:A_2=.04$}
}
\startdata
1 &$ 0.0\pm 0.0242\mathbf{i}$ &$0.0\pm0.9998\mathbf{i}$&$\pm0.2522$&$0.0\pm1.0455\mathbf{i}$&$\pm0.3692$&$0.0\pm1.0926\mathbf{i}$\\

0.8 & $ -1.1614\times10^{-9}\pm0.0254\mathbf{i}$ & $2.5627\times10^{-13}\pm0.9997\mathbf{i}$&$\pm0.2523$& $-3.4918\times10^{-11}\pm1.0455\mathbf{i}$&$\pm 0.3696$& $-7.0611\times10^{-11}\pm1.0927\mathbf{i}$\\

0.6& $-2.8140\times10^{-9}\pm0.0272\mathbf{i} $&$6.9284\times10^{-13}\pm0.9996\mathbf{i}$ &$\pm0.2523$&$-8.4597\times10^{-11}\pm1.0455\mathbf{i}$&$\pm0.36700$&$-1.7126\times10^{-10}\pm1.0928\mathbf{i}$\\

0.4&$-5.5311\times10^{-9}\pm0.0298\mathbf{i}$&$1.5839\times10^{-12}\pm0.9996\mathbf{i}$&$\pm 0.2522 $&$-1.6621\times10^{-10}\pm1.0454\mathbf{i}$&$\pm.3704$&$-3.3695\times10^{-10}\pm1.0929\mathbf{i}$\\

0.2&$-1.1707\times10^{-8}\pm0.0347\mathbf{i}$&$4.3264\times10^{-12}\pm0.9995\mathbf{i}$&$\pm0.2519$&$-3.5125\times10^{-10}\pm1.0454\mathbf{i}$&$\pm0.3707$&$-7.1369\times10^{-10}\pm 1.0930\mathbf{i}$\\
0.0 &Indeterminate&Indeterminate&$0\pm0.0242\mathbf{i}$&Indeterminate&Indeterminate&Indeterminate\\
\enddata
\tablecomments{Table \ref{tbl-7} is presents the roots of characteristic equation(~\ref{eq:cheq}) for $\mu=0.00003$ when P-R effect is included.}
\end{deluxetable}

\end{document}